\documentstyle[12pt]{article}

\newlength{\abstractwidth}
\renewcommand{\thefootnote}{\fnsymbol{footnote}}
\renewcommand{\thanks}[1]{\footnote{#1}} 
\newcommand{\starttext}{ \setcounter{footnote}{0}
\renewcommand{\thefootnote}{\arabic{footnote}}}

\newcommand{\be}{\begin{equation}}
\newcommand{\bea}{\begin{eqnarray}}
\newcommand{\eea}{\end{eqnarray}} \newcommand{\ee}{\end{equation}}
 \newcommand{\<}{\langle}
\renewcommand{\>}{\rangle}
\def\ba{\begin{eqnarray}}
\def\ea{\end{eqnarray}}

\def\cX{{\cal X}}
\def\cL{{\cal L}}
\def\r{\rho}

\def\ra{\rightarrow}

\def\o{\omega}

\def\log{\,{\rm log}\,}

\def\o{\omega}

\def\al{\alpha}
\def\b{\beta}

\def\e{\varepsilon}

\def\l{\lambda}

\def\o{\omega}

\def\r{\rho}

\def\O{\Omega}

\def\ti{\tilde}

\def\C{{\bf C}}

\def\s{\sum}

\def\ddb{{\partial\bar\partial}}

\def\ra{\rightarrow}

\def\us{{\underline s}}

 \def\v{\vskip .1in}

\def\[{{\bf [}}
\def\]{{\bf ]}}

\def\pl{\partial}

\begin{document}
\starttext \baselineskip=15pt \setcounter{footnote}{0}
\newtheorem{theorem}{Theorem}
\newtheorem{lemma}{Lemma}
\newtheorem{definition}{Definition}
\newtheorem{proposition}{Proposition}

\begin{center}
{\Large \bf
REGULARITY OF GEODESIC RAYS
AND MONGE-AMPERE EQUATIONS
\footnote{Work supported in part by the NSF
under grants DMS-07-57372 and
DMS-05-14003.}}
\bigskip

{\large  D.H. Phong and Jacob Sturm} \\

\end{center}

\smallskip
\begin{abstract}

It is shown that the geodesic rays constructed as limits of Bergman geodesics
from a test configuration are always of class $C^{1,\alpha}$, $0<\alpha<1$.
An essential step is to establish that the rays can be extended as solutions of a
Dirichlet problem for a Monge-Amp\`ere equation on
a K\"ahler manifold which is compact.

\end{abstract}

\section{Introduction}
\setcounter{equation}{0}

The purpose of this note is to establish
the $C^{1,\alpha}$ regularity, $0<\al<1$, of
the geodesic rays constructed in \cite{PS07}
from a test configuration by Bergman geodesic approximations.
With the notations given in \S 2 below,
our main result can be stated as follows:

\medskip
\begin{theorem}
\label{maintheorem}
Let $L\to X$ be a positive holomorphic line bundle over a compact
complex manifold $X$. Let $\rho$ be a test configuration
for $L\to X$.
Let $D^\times=\{0<|w|\leq 1\}$ be the punctured unit disk,
and $\pi_X$ the natural projection $X\times D^\times\to X$.
For each metric $h_0$ on $L$
with positive curvature $\o_0\equiv-{i\over 2}\ddb\,\log\,h_0>0$,
let $\Phi(z,w)$ be the $\pi_X^*(\o_0)$-plurisubharmonic
function on $X\times D^\times$
defined by
\bea
\Phi(z,w)={\rm lim}_{k\to\infty}[{\rm sup}_{\ell\geq k}\Phi_\ell(z,w)]^*,
\quad (z,w)\in X\times D^\times,
\eea
where $\Phi_k(z,w)$ are the functions defined by (\ref{bergmank}) below.
Then for any $0<\al<1$,
$\Phi(z,w)$ is a $C^{1,\al}$
generalized solution of the Dirichlet problem
\bea
\label{Dirichlet}
(\pi_X^*(\o_0)+{i\over 2}\ddb\Phi)^{n+1}=0
\ {\rm on}\ X\times D^\times,
\qquad
\Phi(z,w)=0\ {\rm when}\ |w|=1.
\eea
\end{theorem}

\smallskip

The fact that $\Phi(z,w)$ is locally bounded and a solution of the Dirichlet problem
was established in \cite{PS07}, so the new part of the theorem is
the $C^{1,\alpha}$ regularity.
In the case of toric varieties,
the $C^{1,\alpha}$ regularity of geodesic rays
was previously established by Song and Zelditch
\cite{SZ08}, using an explicit analysis of
orthonormal bases for $H^0(X,L^k)$ and the theory
of large deviations. They also pointed out that, already for toric varieties,
geodesic rays from test configurations
can be at best $C^{1,1}$. The interpretation of the completely degenerate Monge-Amp\`ere
equation in (\ref{Dirichlet})
as the equation for geodesics in the space of K\"ahler potentials of class $c_1(L)$
on $X$
is well-known, and due to Donaldson \cite{D99}, Semmes \cite{S}, and Mabuchi \cite{M}.

\smallskip

In \cite{PS09}, $C^{1,\alpha}$ geodesic rays were constructed
in all generality from test configurations by a different approach,
namely viscosity methods for the degenerate complex Monge-Amp\`ere
equation on a compactification $\tilde{\cal X}_D$
of $X\times D^\times$. Thus our theorem can be established by
showing that the above solution, more precisely
$\Phi(z,w)-\Phi_1(z,w)$,
can also be extended to $\tilde{\cal X}_D$
and that such solutions must be unique. For this,
it is essential to show that
$\Phi(z,w)-\Phi_1(z,w)$ is uniformly bounded on $X\times D^\times$.
We accomplish that with the help of a ``lower-triangular'' property
of Donaldson's equivariant imbeddings,
relating $k$-th powers of sections of $H^0(X,L)$ to sections of $H^0(X,L^k)$,
which may be of independent interest (c.f. Lemma \ref{lowertriangular} below).

\smallskip

The uniqueness
follows from a comparison theorem
for Monge-Amp\`ere
equations on K\"ahler manifolds with boundary,
using the
approximation theorems for plurisubharmonic functions obtained recently by
Blocki and Kolodziej \cite{BK} (see also Demailly and Paun \cite{DP}
for other approximation theorems).
It is well-known that such approximation theorems would imply
comparison theorems, by
a straightforward adaptation
to K\"ahler manifolds of the classic comparison theorem of Bedford and Taylor
\cite{BT82} for domains in ${\bf C}^m$. In
fact such a comparison theorem was established in \cite{BK}
for K\"ahler manifolds without boundary. However, the particular version
that we need does not seem available in the literature,
and we have included a brief but complete derivation.

\smallskip

As has been stressed in \cite{PS06}, each test configuration
defines a generalized vector field on the space
of K\"ahler potentials, with the vector at each potential $h_0$
given by the tangent vector $\dot\phi$
to the geodesic at the initial time.
This observation can now be given a precise formulation,
using the measures recently introduced by Berndtsson \cite{B09b}:
for each generalized $C^{1,\alpha}$
geodesic $(-\infty,0]\in t\to\phi(z,t)\equiv \Phi(z,e^t)$, the functional
$\mu_\Phi:C_0^0({\bf R})\ni f
\to \int_X f(\dot\phi)\omega_{\phi(\cdot,t)}^n$
defines a Borel measure on ${\bf R}$ which is independent of $t$.
Taking $t=0$, we can think of this measure as a way of characterizing $\dot\phi(0)$
by its moments.
If $\Phi$ is the geodesic constructed
in Theorem \ref{maintheorem}, the corresponding assignment $h_0\to \mu_\Phi$
can be viewed as a precise realization
of the generalized vector field defined by the test configuration $\rho$.

\smallskip

We note that Theorem \ref{maintheorem} gives the regularity
of the limiting function $\Phi(z,w)$, but it does not provide information
on the precise rate of convergence of $\Phi_k$.
For toric varieties, very precise rates of convergence have been
provided by Song and Zelditch \cite{SZ06, SZ08}. For general manifolds,
in the case of geodesic segments,
the precise rate of $C^0$ convergence
has been obtained a few years ago
by Berndtsson \cite{B09a} with an additional twisting by ${1\over k}K_X$,
and very recently in \cite{B09b} for the $\Phi_k$ themselves.

\smallskip

Finally, we would like to mention that geodesics have been constructed by
Arezzo and Tian \cite{AT}, Chen \cite{C00, C08},
Chen and Tang \cite{CT}, Chen and Sun \cite{CS}, Blocki \cite{B09}
and others in various
geometric situations. For geodesic segments,
the $C^{1,\al}$ regularity has been established by Chen \cite{C00}.
Their construction by Bergman approximations is in \cite{PS06}.
This construction has also been
extended by Rubinstein and Zelditch \cite{RZ} to the construction of harmonic maps in
the space of K\"ahler potentials, in the case of toric varieties.

\section{The extension to a compact K\"ahler manifold}
\setcounter{equation}{0}

In this section, we show how the generalized geodesic rays constructed in \cite{PS07},
originally defined on $X\times \{0<|w|\leq 1\}$,
actually extend as bounded solutions of a complex Monge-Amp\`ere equation
over a compact K\"ahler manifold $\tilde {\cal X}_D\supset X\times \{0<|w|\leq 1\}$.
We begin by introducing the notation and recalling the results of \cite{PS07}.

\subsection{Test configurations}

Let $L\to X$ be a positive line bundle over a compact complex manifold $X$.
A test configuration $\rho$ for $L\to X$ \cite{D02}
is a homomorphism
$\rho:{\bf C}^\times
\to {\rm Aut}({\cal L}\to{\cal X}\to{\bf C})$, where ${\cal L}$
is a ${\bf C}^\times$ equivariant line bundle with ample fibers over a scheme ${\cal X}$,
and $\pi:{\cal X}\to{\bf C}$ is a flat ${\bf C}^\times$ equivariant map of schemes,
with $(\pi^{-1}(1),{\cal L}_{\vert_{\pi^{-1}(1)}})$ isomorphic to $(X,L^r)$ for some fixed $r>0$.
After replacing $L$ by  $L^r$ to a sufficiently high power, we may assume that $r=1$.

\medskip
It is convenient to denote $(\pi^{-1}(w),{\cal L}_{\vert_{\pi^{-1}(w)}})$
by $(X_w,L_w)$. In particular, for each $\tau\not=0$, $\rho(\tau)$
is an isomorphism between $(X_w,L_w)$ and $(X_{\tau w},L_{\tau w})$.

\medskip
The central fiber $(X_0,L_0)$ is fixed under the action of $\rho$. Thus, for each $k$,
$\rho$
induces a one-parameter subgroup of automorphisms
\bea
\rho_k(\tau): H^0(X_0,L_0^k)\to H^0(X_0,L_0^k),
\qquad \tau\in{\bf C}^\times.
\eea
Since $\rho_k(\tau)$ is an algebraic one-parameter subgroup, there is a basis
of $H^0(X_0,L_0^k)$  in which $\r(\tau)$ is represented by a diagonal matrix
with entries
$\tau^{\eta_\al^{(k)}}$,
where $\eta_{\al}^{(k)}$ are integers,
$0\leq\al\leq N_k\equiv {\rm dim}\,H^0(X_0,L_0^k)-1$. Set
\bea
\lambda_\al^{(k)}=\eta_\al^{(k)}-{1\over N_k+1}\sum_{\b=0}^{N_k}\eta_\b^{(k)},
\eea
so that $(\l_\al^{k})$ is the traceless component of $(\eta_\al^{(k)})$.
For a fixed $k$, we shall refer to $\eta_\al^{(k)}$ and $\lambda_\al^{(k)}$
respectively as the weights and the traceless weights
of the test configuration $\rho$.

\smallskip
It is convenient to introduce an $(N_k+1)\times (N_k+1)$ diagonal matrix $B_k$
whose diagonal entries are given by the weights $\eta_\al^{(k)}$.
Such a matrix is determined up to a permutation of the diagonal entries $\eta_\al^{(k)}$,
and we fix one choice once for all. Then the traceless weights
$\lambda_\al^{(k)}$
are the diagonal entries of the matrix $A_k$ defined by
$A_k=B_k-(N_k+1)^{-1}({\rm Tr}\,B_k)I$, and we have
\bea
{\rm Tr}\,B_k=\sum_{\al=0}^{N_k}\eta_\al^{(k)},
\qquad
{\rm Tr}\, A_k=0.
\eea
For sufficiently large $k$, the functions $k(N_k+1)$ and ${\rm Tr}\,B_k$ are
polynomials in $k$ of degree $n+1$, so we have an asymptotic expansion
\bea
{{\rm Tr}\,B_k\over k(N_k+1)}
\equiv F_0+F_1k^{-1}+F_2k^{-2}+\cdots
\eea
The Donaldson-Futaki invariant of
$\r$ is defined to be the coefficient~$F_1$.

\subsection{Equivariant imbeddings of test configurations}

An essential property of test configurations, due to Donaldson
\cite{D05}, is that the entire configuration can be imbedded
equivariantly in ${\bf CP}^{N_k}\times {\bf C}$, in a way which
respects a given $L^2$ metric on $H^0(X,L^k)$. The following formulation
\cite{PS07} is most convenient for our purposes:

\medskip
Let ${\underline s(k)}=\{s_\al^{(k)}(z)\}_{\al=0}^{N_k}$
be a basis for $H^0(X,L^k)$. For all $k$ sufficiently large,
it defines a Kodaira imbedding
\bea
\iota_{\underline s(k)}: X\ni z\to [s_0^{(k)}(z):s_1^{(k)}(z):\cdots:s_{N_k}^{(k)}(z)]\in {\bf CP}^{N_k}
\eea
of $X$ into ${\bf CP}^{N_k}$, with $O(1)$ pulled back to $L^k$.
If $h_0$ is a fixed metric on $L$ with $\o_0\equiv -{i\over 2}\ddb\,\log\,h_0>0$,
then $H^0(X,L^k)$ can be equipped with the $L^2$ metric defined
by the metric $h_0^k$ on sections of $L^k$ and the volume form $\o_0^n/n!$.
For simplicity, we shall refer to this $L^2$ metric on $H^0(X,L^k)$
as just the ``$L^2$ metric defined by $h_0$''.
Of particular importance are then the bases $\underline s(k)$ which are
orthonormal with respect to this $L^2$ metric.

\begin{lemma}
\label{equivariantimbedding}
Let $\rho:{\bf C}^\times\to {\rm Aut}(\cL\to\cX\to{\bf C})$
be a test configuration, and fix a diagonal matrix $B_k$
with the weights of $\rho$ as diagonal entries as defined in \S 2.1.
Fix a metric $h_0$ on $L$ with positive curvature $\o_0$, and corresponding $L^2$
metric on $H^0(X,L^k)$. Then there is an orthonormal basis $\underline s(k)$
of $H^0(X,L^k)=H^0(X_1,L_1^k)$ with respect to the $L^2$ metric defined by
$h_0$ and an imbedding
\bea
I_{\underline s}:
({\cal L}\to{\cal X}\to{\bf C})
\to
(O(1)\times{\bf C}\to{\bf CP}^{N_k}\times {\bf C}\to {\bf C})
\eea
satisfying

(1) $I_{{\underline s}(k)}\vert_X=\iota_{\us(k)}$

(2) $I_{\us(k)}$ intertwines $\rho(\tau)$ and $B_k$,
\bea
I_{\us(k)}(\rho(\tau)\ell_w)
=
(\tau^{B_k}I_{\us(k)}(\ell_w),
\tau w),
\qquad \ell_w\in L_w,\ \tau\in{\bf C}^\times.
\eea
\end{lemma}

Let $E_k=\pi_*(\cL^k)$ be the direct images of the bundles $\cL^k$.
Thus $E_k\to{\bf C}$
is a vector bundle over ${\bf C}$ of rank $N_k+1$, and its sections $S(w)$
are holomorphic sections of $L_w$ for each $w\in{\bf C}$. The action
of ${\bf C}^\times$ on the sections $S$ is given by
\bea
S^\tau(w)=\rho(\tau)^{-1}S(w\tau).
\eea
Then a third key statement in the equivariant imbedding lemma
is:

{\it (3) The functions
\bea
S_\al(w)\equiv w^{\eta_\al^{(k)}}\rho(w)\,s_\al,
\qquad w\in{\bf C}^\times
\eea
extend to a basis for the free ${\bf C}[w]$ module of all sections
of $E_k\to{\bf C}$ and they have the property: $S_\al(1)=s_\al$. This extension still satisfies
the relation}
\bea
\rho(\tau)^{-1}S_\al(w)=\tau^{\eta_\al^{(k)}}S_\al(w),
\qquad w\in{\bf C}.
\eea

\subsection{The construction of geodesics}

We come now to the construction of geodesics by Bergman approximations.
Let
$\rho:{\bf C}^\times
\to {\rm Aut}({\cal L}\to{\cal X}
\to{\bf C})$
be a test configuration for $L\to X$, and fix a metric $h_0$ on $L$
with positive curvature $\o_0$. Let
$\us(k)=\{s_\al^{(k)}(z)\}$ be an orthonormal basis for $H^0(X,L^k)$ with
respect to the $L^2$ metric defined by $h_0$ as in Lemma \ref{equivariantimbedding}.
Let
\bea
D^\times=\{w\in {\bf C};\ 0<|w|\leq 1\}
\eea
be the punctured disk. Define the functions $\Phi_k(z,w)$ by
\bea
\label{bergmank}
\Phi_k(z,w)
=
{1\over k}\log\,\sum_{\al=0}^{N_k}
|w|^{2\eta_\al^{(k)}}|s_\al^{(k)}(z)|_{h_0^k}^2-{n\over k}\log\,k,
\qquad (z,w)\in X\times D^\times.
\eea
and $\Phi(z,w)$ by
\bea
\label{bergman}
\Phi(z,w)={\rm lim}_{k\to\infty}[{\rm sup}_{\ell\geq k}\Phi_\ell(z,w)]^*
\eea
where $\eta_\al^{(k)}$
are the weights of the test configuration $\rho$,
${}^*$ denotes the upper semi-continuous envelope,
i.e. $f^*(z)={\rm lim}_{\epsilon\to 0}{\rm sup}_{|w-z|<\epsilon}f(w)$,
and $|s_\al(z)|_{h_0^k}^2\equiv s_\al(z)\overline{s_\al(z)}h_0(z)^k$
denotes the norm-squared of $s_\al(z)$ with respect to the metric $h_0^k$.
Then it is shown in \cite{PS07}
\footnote{Actually, in \cite{PS07}, the weights $\eta_\al^{(k)}$ in
the definition of $\Phi_k(z,w)$ were replaced by the traceless weights
$\lambda_\al^{(k)}$. If we denote by $\Phi_k^\#(z,w)$
the functions obtained in this manner with the traceless weights,
then we have
\bea
\Phi_k(z,w)=\Phi_k^\#(z,w)+{{\rm Tr}\,B_k\over k(N_k+1)}\log\,|w|^2.
\eea
It follows that the complex Hessians of $\Phi_k(z,w)$ and $\Phi_k^\#(z,w)$
are identical. However, the behaviors near $|w|=0$ of $\Phi_k(z,w)$
and $\Phi_k^\#(z,w)$ are different, and for our purposes, it is important
to work with $\Phi_k(z,w)$.}
that $\Phi(z,w)$ is a generalized geodesic ray in the sense that

\smallskip

(a) $\pi_X^*(\o_0)+{i\over 2}\ddb\Phi\geq 0$ on $X\times D^\times$,
where $\pi_X$ is the projection $X\times D^\times\to X$ on the first factor;

(b) For each finite $T>0$, we have
\bea
\label{Tbound}
{\rm sup}_k|\Phi_k(z,w)|,\
|\Phi(z,w)|\leq C_T
\qquad
{\rm for}\ (z,w)\in X\times \{e^{-T}<|w|\leq 1\}
\eea
with $C_T$ a constant independent of $z,w$ and $k$,
but possibly depending on $T$;

(c) $\Phi(z,w)$ is continuous when $|w|=1$, and is a solution
in the sense of pluripotential theory of the following Dirichlet
problem for the completely degenerate Monge-Amp\`ere equation
\bea
\label{dirichlet}
(\pi_X^*(\o_0)+{i\over 2}\ddb\Phi)^{n+1}=0
\ {\rm on}\ X\times D^\times,
\qquad
\Phi(z,w)=0 \ {\rm when}\ |w|=1.
\eea

\smallskip
The geodesic $\Phi(z,w)$ is non-constant if the test configuration is non-trivial,
that is, not holomorphically equivalent to a product test configuration.
We note that in the boundary value problem (\ref{dirichlet}),
the behavior of $\Phi(z,w)$ near $w=0$ is not specifically assigned.

\subsection{Formulation in terms of equivariant imbeddings}

We come now to the main task in this chapter, which is to identify
the solution (\ref{dirichlet}) with the restriction
to $X\times D^\times$ of the solution of
a standard Dirichlet problem on a compact
K\"ahler manifold $\tilde{\cal X}_D$ with boundary.

\medskip
Let $\pi:\cX\ra\C$ be the projection map, and $D=\{w\in\C: |w|\leq 1\}$.
Let $\cX_D=\pi^{-1}(D)$, ${\cal X}_D^\times=\pi^{-1}(D^\times)$.
The space ${\cal X}_D^\times$
is isomorphic to $X\times D^\times$ under the correspondence
\bea
\label{trivialization}
X\times D^\times\ni (z,w)\to \rho(w)(z)
\in  X_w,
\eea
where $z\in X$ is viewed as a point in $X_1$. This correspondence lifts to
a correspondence between $L\times D^\times$ and the restriction
${\cal L}_D^\times$ of ${\cal L}$ over ${\cal X}_D^\times$.
\v
Let $p:\ti\cX\ra\cX\to {\bf C}$ be an $S^1$ equivariant smooth
resolution and $\ti\cL=p^*\cL$.
The first step is to show that the functions $\Phi_k(z,w)-\Phi_1(z,w)$
of (\ref{bergmank}),
which are defined on $X\times D^\times$,  may be
extended to
plurisubharmonic functions on all of $\ti\cX_D=p^{-1}(\cX_D)$.

\smallskip
Let us fix a metric $h_0$ on $L$ with positive curvature $\o_0$.
Let $\us(k)$ be the orthonormal basis for $H^0(X,L^k)$ with respect to
$h_0$ provided by Lemma
\ref{equivariantimbedding}, and let $I_{\us(k)}$ be a corresponding
equivariant imbedding of the test configuration. Let $\Phi_k(z,w)$ be defined by (\ref{bergmank}).
Define a closed $(1,1)$-form $\Omega_k$ on $\ti\cX_D$ by
\bea
\label{Phik}
\Omega_k
=
{1\over k} (I_{\us(k)}\circ p)^*\o_{FS}
\eea
where $\o_{FS}$ is the Fubini-Study metric on ${\bf CP}^{N_k}$. Define as well
a hermitian metric $H_k$ on $\ti\cL$
by $H_k = (I_{\us(k)}\circ p)^*
(h_{\rm FS})^{1/k}$, where $h_{\rm FS}$
is the Fubini-Study metric on the hyperplane bundle $O(1)$ over
${\bf CP}^{N_k}$.
Thus $\Omega_k$ is the
curvature of $H_k$. The restriction of
$\O_k$ to ${\cal X}_D^\times$ can be readily worked out explicitly in terms of the coordinates
$(z,w)$. Using the intertwining property of the equivariant imbedding,
\bea
X\times D^\times\ni(z,w)
\to \rho(w)z
\to
&
I_{\us(k)}(\rho(w)z)=w^{B_k}I_{\us(k)}(z)
=
(w^{B_k}\iota_{\us(k)}(z),w)
\eea
we find that $I_{\us(k)}$ is given by
\bea
I_{\us(k)}:
X\times D^\times
\ni
(z,w)
\to
([w^{\eta_0^{(k)}}s_0^{(k)}(z):w^{\eta_1^{(k)}}s_1^{(k)}(z):
\cdots:
w^{\eta_{N_k}^{(k)}}s_{{N_k}}^{(k)}(z)],w).
\eea
Since the Fubini-Study metric $h_{\rm FS}$ on $O(1)$
at $[s_0:s_1:\cdots:s_{N_k}]\in{\bf CP}^{N_k}$ is given by
$h_{\rm FS}=(|s_0|^2+\cdots+|s_{N_k}|^2)^{-1}$, we obtain the following expression for
$\Omega_k$,
\bea
\Omega_k|_{X\times D^\times}=
{1\over k}\dot {i\over 2}\ddb\,\log\,\sum_{\al=0}^{N_k}
|w|^{2\eta_\al^{(k)}}|s_\al^{(k)}(z)|^2.
\eea
Recalling that the norm with respect to $h_0^k$
of a section $s(z)$ of $L^k$ is given
by $|s(z)|_{h_0^k}^2=|s(z)|^2 h_0^k$, we
find the following key relation between the $(1,1)$-forms
$\Omega_k$ and
the potentials $\Phi_k(z,w)$ defined earlier
in (\ref{bergmank}),
\bea
\label{OmegaPhik}
\Omega_k|_{X\times D^\times}
\ =\
\pi_X^*(\o_0)+{i\over 2}\ddb\Phi_k(z,w).
\eea

\subsection{The extension of $\Psi_k$ to the total space $\ti\cX_D$}

The relation (\ref{OmegaPhik}) that we have just obtained shows
that the form $\pi_X^*(\o_0)+{i\over 2}\ddb\Phi_k(z,w)$,
defined originally on ${\cal X}_D^\times$, admits the natural extension
$\Omega_k$ to the whole of $\ti{\cal X}_D$.

\smallskip
Since the form $\pi_X^*(\o_0)$ does not extend by itself to $\tilde{\cal X}$,
we re-write $\Omega_k$ as
\bea
\Omega_k=\Omega_1
+
{i\over 2}\ddb\,(\Phi_k-\Phi_1)
\equiv
\Omega_1+{i\over 2}\ddb\,\Psi_k.
\eea
The function $\Psi_k=\Phi_k-\Phi_1$ has a simple interpretation that shows that it
extends as a smooth function to the whole of $\ti\cX_D$: as we saw earlier
in \S 2.2,
under the maps $I_{\us(k)}$ and $I_{\us(1)}$
of the test configuration $\rho$,
the Fubini-Study metric $h_{\rm FS}$ pulls back respectively
to $H_k^k=(\sum_\al |w|^{2\eta_\al^{(k)}}|s_\al^{(k)}(z)|^2)^{-1}$
and $H_1=(\sum_\al|w|^{2\eta_\al^{(1)}}|s_\al^{(1)}(z)|^2)^{-1}$ on $L\times D^\times$.
Thus we may
write
\bea
\Psi_k=\log {H_1
\over
H_k}-{n\over k}\log\,k.
\eea
The right hand side is a well-defined, smooth scalar function over the whole
of $\ti\cX_D$, since it is the logarithm of the ratio of two smooth metrics
on the same line bundle $\ti\cL\to\ti\cX_D$.

Since $\Omega_k$ is non-negative as the pull-back of a non-negative form, the function
$\Psi_k$ is $\Omega_1$-plurisubharmonic
\footnote{In general, given a non-negative smooth, closed $(1,1)$-form $\Omega$
on a complex manifold $X$, we say that a scalar function $\Phi$
is $\Omega$-plurisubharmonic if $f_\al+\Phi$
is plurisubharmonic on $U_\al$ for each $\al$,
if $X=\cup_\al U_\al$ is a covering of $X$ by coordinate charts $U_\al$
with $\Omega={i\over 2}\ddb f_\al$ on $U_\al$.}.
We also define
\bea
\label{Psi}
\Psi={\rm lim}_{k\to\infty}[{\rm sup}_{\ell\geq k}\Psi_\ell]^*
\eea
which is an extension of $\Phi-\Phi_1$ to $\ti\cX_D$.

\subsection{Uniform estimates for $\Psi_k$}

Recall that in \cite{PS07}, as quoted in (\ref{Tbound}) above,
we only have bounds for the functions $\Phi_k(z,w)$
when $|w|>e^{-T}$, for some fixed finite $T>0$. Since the
function $\Psi_k$ extends to a smooth function on $\ti\cX_D$,
it follows that it is bounded on $\ti\cX_D$. However, the bound may
a priori depend on $k$. The most important step
in the extension to $\ti\cX_D$ is to show that this bound can
actually be made uniform in $k$.

\medskip
We carry this out with several lemmas. The first is the following
essential ``lower triangular lemma'':

\begin{lemma}
\label{lowertriangular}
Fix a test configuration $\rho$, and a metric $h_0$ on $L$ with
positive curvature $\o_0$. For each $k$, let $\us(k)=\{s_\al^{(k)}\}_{\al=0}^{N_k}$
be an orthonormal basis for $H^0(X,L^k)$ as in Lemma \ref{equivariantimbedding}.
Then for any $s_\b^{(1)}$ in $\us(1)$, we can write
\bea
\label{lowertriangularinequality}
(s_\b^{(1)})^k
=
\sum_{\l_\al^{(k)}\leq k\l_\b^{(1)}}a_{\b\al}\,s_\al^{(k)}
\eea
where $a_{\b\al}\in\C$  and the subindex indicates the range of
indices $\al$ which are allowed. Furthermore, the coefficients $a_{\b\al}$
satisfy the bound
\bea
|a_{\b\al}|
\leq V^{1\over 2}M^k
\eea
where we have set $M={\rm sup}_{0\leq \b\leq N_1}{\rm sup}_X|s_\b^{(1)}|_{h_0}$
and $V=\int_X \o_0^n$.
\end{lemma}

\noindent
{\it Proof of Lemma \ref{lowertriangular}}:
For each $k$, let $E_k=\pi_*(\cL^k)\to{\bf C}$,
and let $S_0(w),\cdots,S_{N_k}(w)$ be a basis for the free ${\bf C}[w]$
module of sections of $E_k\to{\bf C}$,
as provided in Lemma \ref{equivariantimbedding}.
Now let $S_\beta$ be an element of this basis for $E_1\ra \C$, and some $\beta$ with
$0\leq\beta\leq N_1$. Then $\r(\tau)^{-1}S_\b(w\tau)=\tau^{\eta_\b^{(1)}}S_\b(w)$ which implies
\be
\r(\tau)^{-1}S_\b^k(w\tau)=\tau^{k\eta_\b^{(1)}}S^k_\b(w)
\ee
On the other hand, $S_\b^k$ is a section of $E_k$ so we may write
\be
\label{free} S_\b^k(w)\ = \ \s_{\al=0}^{N_k}\, a_\al(w)S_\al(w)
\ee
for certain uniquely defined polynomials $a_\al(w)\in\C[w]$. Applying the $\C^\times$ action to
both sides of (\ref{free}) we obtain
\be
\s_{\al=0}^{N_k}\tau^{k\eta_\b^{(1)}}a_\al(w)S_\al(w)\ = \ \tau^{k\eta_\b^{(1)}}S_\b^k(w)\ = \
\r(\tau)^{-1}S_\b^k(w\tau)\  = \  \s_{\al=0}^{N_k} a_\al(w\tau)\tau^{\eta_\al^{(k)}}S_\al(w)
\ee
Comparing coefficients we obtain
\be
\label{coef}
\tau^{k\eta_\b^{(1)}}a_\al(w)\ = \ a_\al(w\tau)\tau^{\eta_\al^{(k)}}
\ee
Setting $w=1$ we see that $a_\al(\tau)=a_{\b\al}\tau^{r_\al}$ for some integer $r_\al$
and some $a_{\b\al}\in\C$.
But $a_\al(w)$ is
a polynomial. Thus $r_\al\geq 0$ and $a_\al(w)=a_{\b\al} w^{r_\al}$ for all $w\in\C$. The
equation (\ref{coef})
implies that if $a_{\b\al}\not=0$ we have $k\eta_\b^{(1)}= r_\al+\eta_\al^{(k)}$ and thus $\eta_\al^{(k)}\leq k\eta_\b^{(1)}$.
Evaluating (\ref{free}) at $w=1$ we obtain the first part of the lemma.

\smallskip
Finally, the orthonormality of the sections $s_\al^{(k)}$ implies
\bea
|a_{\b\al}|
=
|\int \<(s_\b^{(1)})^k,s_\al^{(k)}\>_{h_0^k}\o_0^n|
\leq
\int |s_\b^{(1)}|_{h_0}^k\cdot |s_\al^{(k)}|_{h_0^k}\o_0^n
\leq M^k\,V^{1\over 2}
\eea
The lemma is proved.
\v

\noindent
{\it Remark}: It may happen that
$\eta_\al^{(k)}<k\eta_\b^{(1)}$ for all $\al$ with $a_{\b\al}\not=0$, that is, it may
happen that $a_\al(w)$ vanishes at $w=0$ for all $\al$. This would mean
that $S_\b(0)$ is a non-zero section of $H^0(X_0,L_0)$ but that $S_\b^k(0) = 0\in H^0(X_0,L_0^k)$,
in other words,
the section $S_\b(0)$ is nilpotent (which is possible if $X_0$ is a non-reduced scheme, that
is, if $X_0$ has nilpotent elements in its structure sheaf).
\v

Next, we also need

\begin{lemma}
\label{kaehler}
The complex manifold $\ti\cX_D$ always admits a K\"ahler metric.
\end{lemma}

\medskip
\noindent
This lemma was proved in \cite{PS07a}. In fact, it is proved there that there
exists a line bundle ${\cal M}$ on $\cX_D$
which is trivial on $\cX^\times$, and such that $\cL^m\otimes {\cal M}$
is positive for some fixed positive power $m$. The desired K\"ahler metric on
$\ti\cX_D$ can then be taken to be the ratio of the curvature of $\cL^m\otimes{\cal M}$
by $m$. Q.E.D.

\begin{lemma}
\label{infinitybound}
There exists a finite constant $C$ so that
\bea
\label{infinityboundinequality}
{\rm sup}_{k\geq 1}{\rm sup}_{\tilde\cX_D}|\Psi_k|\leq C<\infty.
\eea
In particular,
\bea
{\rm sup}_{\tilde\cX_D}|\Psi|\leq C<\infty.
\eea
\end{lemma}

\noindent
{\it Proof of Lemma \ref{infinitybound}}:
Let $H$ be a K\"ahler metric on $\tilde\cX_D$, which exists by Lemma \ref{kaehler}.
Since $\Psi_k$ is $\Omega_1$-plurisubharmonic, it follows that
\bea
\Delta_H\Psi_k\geq -C_1,
\eea
where $\Delta_H$ is the Laplacian with respect to $H$, and $C_1$
is an upper bound for the trace of $\Omega_1$ with respect to
the metric $H$. On the other hand, $\Psi_k\vert_{\pl\tilde\cX_D}\to -\Phi_1$
uniformly as $k\to\infty$, and thus $\Psi_k\vert_{\pl\tilde\cX_D}\leq C_2$.
Let $u$ be the smooth function on $\tilde\cX_D$ which is the solution of the
Dirichlet problem
\bea
\Delta_Hu=-C_1 \
{\rm on}\ \tilde\cX_D,
\qquad
u=C_2\
{\rm on}\ \pl\tilde\cX_D.
\eea
By the maximum principle, we have $\Psi_k\leq u$ for all $k$,
and this gives the upper bound.

\smallskip

To establish the lower bound,
it suffices to prove that
\bea
\Psi_k \geq -C
\
{\rm on}
\ \cX_D^\times
\eea
where $C$ is a constant independent of
$k$, since each function $\Psi_k$ is smooth on $\tilde\cX_D$.
On $\cX^\times$, we can use the explicit expressions
for $X\times D^\times$ and write
\bea
\Psi_k
=
\log {(\sum_{\al=0}^{N_k}|w|^{2\eta_\al^{(k)}}|s_\al^{(k)}|_{h_0^k}^2)^{1\over k}
\over
\sum_{\b=0}^{N_1}|w|^{2\eta_\b^{(1)}}|s_\b^{(1)}|_{h_0}^2}
-
{n\over k}\log k.
\eea
Now fix $w$ with $0<|w|\leq 1$, fix
$z\in X$, and choose $\b_0$ so that
\bea
|w|^{2\eta_{\b_0}^{(1)}}|s_{\b_0}^{(1)}(z)|_{h_0}^2
=
{\rm sup}_{0\leq\b\leq N_1}
|w|^{2\eta_{\b}^{(1)}}|s_{\b}^{(1)}(z)|_{h_0}^2.
\eea
In view of Lemma \ref{lowertriangular}, we can write
\bea
|(s_{\b_0}^{(1)})^k|_{h_0^k}
\leq M^k V^{1\over 2}\sum_{k\eta_{\b_0}^{(1)}\geq \eta_{\al}^{(k)}}
|s_\al^{(k)}|_{h_0^k}.
\eea
Since $|w|\leq 1$, we have then
\bea
|w|^{2k\eta_{\b_0}^{(1)}}|s_{\b_0}^k(z)|_{h_0^k}^2
&\leq&
M^{2k}V (\sum_{k\eta_{\b_0}^{(1)}\geq \eta_{\al}^{(k)}}|w|^{\eta_\al^{(k)}}|s_\al^{(k)}(z)|_{h_0^k})^2
\nonumber\\
&\leq&
M^{2k}V
(N_k+1)
\sum_{k\eta_{\b_0}^{(1)}\geq \eta_{\al}^{(k)}}|w|^{2\eta_\al^{(k)}}|s_\al^{(k)}(z)|_{h_0^k}^2.
\eea
Returning to $\Psi_k$, we can now write
\bea
\Psi_k(z,w)
&\geq&
\log {(\sum_{\al=0}^{N_k}|w|^{2\eta_{\al}^{(k)}}|s_\al^{(k)}(z)|_{h_0^k}^2)^{1\over k}
\over
(N_1+1)|w|^{2\eta_{\b_0}^{(1)}}|s_{\b_0}^{(1)}(z)|_{h_0}^2}
-
{n\over k}\log k
\nonumber\\
&\geq&
-{1\over k}\log (V(N_k+1))
-2 \log M-{n\over k}\log k-\log (N_1+1)
\eea
in view of the preceding inequality. This establishes
Lemma \ref{infinitybound} since $N_k\leq C\,k^n$.

\subsection{The Monge-Amp\`ere equation on the whole of $\tilde\cX_D$}

With the uniform estimates provided by Lemma \ref{infinitybound},
it follows readily that the function $\Psi$ defined by
(\ref{Psi}) is a bounded, $\Omega_1$-plurisubharmonic function
on $\tilde\cX_D$.
Since it satisfies a completely degenerate Monge-Amp\`ere equation on $\cX_D^\times$,
and since the singular set $X_0$ is an analytic subvariety,
it follows from general pluripotential theory that it satisfies the same
completely degenerate equation on $\tilde\cX_D$. We give now a direct proof
of this fact, since we already have at hand all
the necessary ingredients. It suffices to observe that $\Psi_k$
satisfies the following properties:

\begin{lemma}
\label{Psikproperties}
The functions $\Psi_k$ satisfy

(a) ${\rm sup}_k\,{\rm sup}_{\ti\cX_D}|\Psi_k|\leq C<\infty$;

(b) $\int_{\ti\cX_D}(\Omega_{1}+{i\over 2}\ddb\Psi_k)^{n+1}\leq C{1\over k}$;

(c) Let $T$ be the vector field $T={\pl\over\pl t}$
defined in a neighborhood of the boundary $|w|=1$ on $\tilde\cX_D$, where $t=\log |w|$.
Then ${\rm sup}_U|T\Psi_k|\leq C$, where $C$ is a constant, and $U$
is a neighborhood of the boundary $|w|=1$, independent of $k$.

(d) ${\rm sup}_{\pl\tilde\cX_D}|\Psi_k+\Phi_1|\leq a_k$, with $a_k$ decreasing to
$0$ and $\sum_{k=1}^\infty a_k<\infty$.

\end{lemma}

\noindent
{\it Proof of Lemma \ref{Psikproperties}}: Part (a) is just the statement
of Lemma \ref{infinitybound}. Part (b) follows from the fact that
the form $\Omega_{1}+{i\over 2}\ddb\Psi_k$ is smooth on $\cX_D$,
and that its Monge-Amp\`ere mass on $\cX_D^\times$ coincides with
the Monge-Amp\`ere mass of $\O_k=(\pi_X^*(\o_0)+{i\over 2}\ddb\Phi_k)$
on $X\times D^\times$. As we already observed in the footnote 1,
$\Phi_k$ and $\Phi_k^\#$ have the same complex Hessian. So the desired
estimate follows from the analogous estimate for the Monge-Amp\`ere mass
of $(\pi_X^*(\o_0)+{i\over 2}\ddb\Phi_k^\#)$ established in Lemma 4.3
of \cite{PS07}. Part (c) follows from the bound $|\eta_\al^{(k)}|\leq C\,k$,
established in Lemma 3.1 of \cite{PS07}. Finally, Part (d), with
$a_k=C\,k^{-2}$, follows
from the Tian-Yau-Zelditch theorem \cite{T90,Y93,Z} (see also Catlin \cite{Ca}
and Lu \cite{L})
as shown in the case of geodesic segments in \cite{PS06}.
Q.E.D.

\medskip
We can now formulate the main theorem of this chapter:

\begin{theorem}
\label{extensiontheorem}
Let $L\to X$ be a positive line bundle over a compact complex manifold,
let $\rho$ be a test configuration, and let $h_0$ be a metric on $L$ with
positive curvature $\o_0$. Let $\tilde\cX$
be an $S^1$ invariant resolution $p:\tilde\cX\to\cX\to{\bf C}$ of $\cX$,
and $\tilde\cX_D=(\pi\circ p)^{-1}(D)$.
Let $\Phi_k, \Phi$ be defined as in (\ref{bergmank}) and
(\ref{bergman}). Set
\bea
\Psi\ =\ \Phi-\Phi_1\
{\rm on}\
X\times D^\times.
\eea
Then the function $\Psi$ extends to a bounded, $\Omega_1$-plurisubharmonic function
on $\tilde\cX_D$, which is a generalized solution of the following
Dirichlet problem on $\ti\cX_D$,
\bea
\label{dirichletcX}
(\Omega_{1}+{i\over 2}\ddb\Psi)^{n+1}=0
\
{\rm on}\ \tilde\cX_D,
\qquad
\Psi=-\Phi_1\
{\rm on}
\ \pl\tilde\cX_D.
\eea
Here $\Omega_1$ is the pull-back to
$\tilde\cX_D$ of the Fubini-Study metric by $I_{\us(1)}\circ p$.
\end{theorem}

\subsection{Positivity of the background form away from the central
fiber}

The equation (\ref{dirichletcX}) provides an extension of the degenerate complex
Monge-Amp\`ere equation to the compact manifold with boundary $\tilde{\cal X}_D$.
It is however written with respect to a background $(1,1)$-form
$\Omega_1$ which may be degenerate.
In preparation for applications of uniqueness theorems for the
complex Monge-Amp\`ere equation, we rewrite it now with a background $(1,1)$-form
which is non-negative everywhere, and strictly positive away from
the central fiber $p^{-1}(X_0)$.

\smallskip
For this, we make use of Lemma 1 of \cite{PS09}, which asserts the existence
of a $S^1$ invariant metric $H_0$ on $\tilde\cL$ with the following properties:
\bea
\label{H0}
&&
\Omega_0\equiv-{i\over 2}\ddb\log\,H_0
\geq 0
\quad
{\rm on}
\
\tilde\cX_D,
\quad
\Omega_0 >0
\ \
{\rm on}\
\cX_D^\times
\nonumber\\
&&
H_0\vert_{\pl\tilde\cX_D}=h_0.
\eea

Let $\Psi_0$ be defined by
\bea
\label{Psi0}
\Psi_0=\log {H_0\over (I_{\us(1)}\circ p)^*(h_{\rm FS})}=\log {H_0\over H_1},
\eea
which is a smooth function on $\tilde\cX_D$, since it is the logarithm of the ratio of two metrics
on the same line bundle $\tilde\cL$. Restricted to $\pl\tilde\cX_D$, it is given by
\bea
\Psi_0\vert_{\pl\tilde\cX_D}
=
\log {h_0\over (\sum_{\al=0}^{N_1}|s_\al^{(1)}|^2)^{-1}}
=
\log \sum_{\al=0}^{N_1}|s_\al^{(1)}|_{h_0}^2
=
\Phi_1\vert_{\pl\tilde\cX_D}.
\eea
Let $\Psi$ be the solution on $\tilde\cX_D$ of the completely degenerate Monge-Amp\`ere equation
with background form $\Omega_1$ as given in Theorem \ref{extensiontheorem}.
Define the function $\hat\Phi$ on $\tilde\cX_D$ by
\bea
\label{hatPsi}
\hat\Phi=\Psi+\Psi_0.
\eea
Clearly $(\Omega_0+{i\over 2}\ddb \hat\Phi)^{n+1}=0$ on $\tilde\cX_D$. Furthermore,
restricted to the boundary $\pl\tilde\cX_D$, we have
\bea
\hat\Phi\vert_{\pl\tilde\cX_D}=
\Psi\vert_{\pl\tilde\cX_D}+\Psi_0\vert_{\pl\tilde\cX_D}
=
-\Phi_1\vert_{\pl\tilde\cX_D}+\Phi_1\vert_{\pl\tilde\cX_D}=0.
\eea
In summary, we have obtained the following alternative formulation of
Theorem \ref{extensiontheorem}:

\begin{theorem}
\label{MAOmega0theorem}
Let the setting be the same as in Theorem \ref{extensiontheorem},
and let $H_0$ be a metric on $\cL$ as in (\ref{H0}),
$\Psi_0$ be defined as in (\ref{Psi0}), and $\hat\Psi\equiv \Phi-\Phi_1+\Psi_0$.
Then the function $\hat\Psi$ is a bounded, $\Omega_0$-plurisubharmonic
generalized solution
of the following Dirichlet problem,
\bea
\label{MAOmega0}
(\Omega_0+{i\over 2}\ddb\hat\Phi)^{n+1}=0
\
{\rm on}
\ \tilde\cX_D,
\qquad
\hat\Phi\vert_{\pl\tilde\cX_D}=0.
\eea
\end{theorem}

\section{The comparison principle on K\"ahler manifolds}
\setcounter{equation}{0}

We derive now the uniqueness theorem that we need. We note that
there has been considerable progress recently on uniqueness theorems
for the complex Monge-Amp\`ere equation, and in particular for
certain broad classes of possibly unbounded solutions (see e.g. Blocki \cite{B03},
Blocki and Kolodziej \cite{BK}, Dinew \cite{D}, and references therein). However, there does not
appear to be a version that would apply directly to our situation,
namely to the Dirichlet problem on K\"ahler manifolds with boundary, for
$\Omega_0$-plurisubharmonic functions where the closed $(1,1)$-form
$\Omega_0$ is non-negative, but may be degenerate. We provide such a
version below, just by following the original arguments of
Bedford and Taylor \cite{BT82} in ${\bf C}^m$.

\subsection{The comparison principle}

\begin{theorem}
Let $(M,\Omega)$ be a compact K\"ahler manifold with smooth boundary
$\partial M$ and dimension $m$, and let $\Omega_0$ be a smooth, non-negative, closed
$(1,1)$-form. Then we have
\bea
\int_{\{u<v\}}(\Omega_0+{i\over 2}\ddb v)^m
\leq
\int_{\{u<v\}}(\Omega_0+{i\over 2}\ddb u)^m.
\eea
for all $u,v$ in $L^\infty$,
$\Omega_0$-plurisubharmonic,
and satisfying ${\rm liminf}_{z\to\partial M}(u(z)-v(z))\geq 0$.
\end{theorem}

\medskip

We adapt the original proof of Bedford-Taylor \cite{BT82} to our setting.
The main steps are as follows. The first step is a version of the theorem,
in the special case of smooth data:

\begin{lemma}
\label{BT76}
Let $u,v\in C^\infty(\bar M)$ be $\Omega_0$-plurisubharmonic functions satisfying
$u(z)-v(z)\geq 0$ for $z\in\partial M$,
and $\Omega_0$ a smooth, closed, non-negative $(1,1)$-form.
Assume that $\{u<v\}$ has smooth
boundary. Then
\bea
\int_{\{u<v\}}(\Omega_0+{i\over 2}\ddb v)^m
\leq
\int_{\{u<v\}}(\Omega_0+{i\over 2}\ddb u)^m.
\eea
\end{lemma}

The proof is identical to \cite{BT76}, Proposition 4.1.
Next, we need a notion of capacity adapted to
K\"ahler manifolds, as e.g. in \cite{PS06}:

\smallskip
Let $M'$ be any open subset
with compact closure in $M$,
and let $M'\subset\cup_{\al=1}^N U_\al$ be a finite cover of $M'$ by
a fixed system of coordinate neighborhoods $U_\al$.
We say that $E\subset M'$ has capacity $c(E,N)<\e$ if we can write
$E=\cup_{\al=1}^N E_\al$, with $E_\al\subset U_\al$ Borel subsets
and $\sum_{\al=1}^N c(E_\al, U_\al)<\e$, where
$c(A,B)$ is the capacity for subsets of ${\bf C}^m$,
\bea
c(A,B)
=
{\rm sup}\,\{\int_A ({i\over 2}\ddb v)^m;\ v\ {\rm plurisubharmonic},
\ 0\leq v\leq 1\}.
\eea
We say that $c(E,M)=0$ if $c(E,M)<\e$ for every $\e>0$.
With this definition, it is easy to extend the quasi-continuity
theorem of Bedford-Taylor \cite{BT82} to K\"ahler manifolds:
any $\Omega_0$-plurisubharmonic function $u$ on $M$
is ``quasi-continuous'', i.e.,
for any open subset $M'$ with compact closure and
any $\e>0$, there is an open set $G\subset M'$
so that $c(G,M')<\e$ and $u$ is continuous on $M'\setminus G$.
With this notion of capacity, we can axiomatize
the limiting processes in \cite{BT82}:

\begin{lemma}
\label{BT82-1}
Let $M'$ be a fixed open subset with
compact closure of the complex manifold $M$.
Assume that $u,v,u_k,v_j$
are Borel measurable functions on an open neighborhood of $\bar M'$,
and $d\mu_k, d\nu_j, d\mu, d\nu$ are non-negative Borel measures
on $M'$ with the following properties:

\smallskip

(a) $u,v$ are upper semi-continuous and quasi-continuous;

(b) $u_j,v_j\in C^\infty(M')$,
$u_j$, $v_j$ decrease to $u$ and $v$ respectively. Furthermore,
there exists a $\delta>0$ and a neighborhood $W$
of $M\setminus M'$ so that,
for any $k$, there exists $J_k$ satisfying
\bea
u_k\geq v_j+\delta
\ {\rm on}\ W,
\ {\rm for}\ j\geq J_k.
\eea

(c) $d\mu_k\to d\mu$, $d\nu_j\to d\nu$
weakly as $k$ and $j$ tend to $\infty$;

(d) The measures $d\mu_k$, $d\nu_k$, $d\mu$, $d\nu$ are uniformly bounded
with respect to capacity, in the following sense:
there exists a constant $C$ such that for any $\e$, $0<\e<1$
and any Borel subset $G\subset M'$ with
$c(G,M')<\e$, 
\bea
\int_{G}(d\mu_k+d\nu_k+d\mu+d\nu)\leq C\,\e
\qquad {\rm for\ all}\ k.
\eea

(e) For each $k$, there exists $J_k$ so that
\bea
\label{comparisonineq1}
\int_{\{u_k<v_j\}}
d\nu_j
\leq
\int_{\{u_k<v_j\}}
d\mu_k
\quad {\rm for}\ j\geq J_k.
\eea
Then we can conclude that
\bea
\label{comparisonineq2}
\int_{\{u<v\}}d\nu
\leq
\int_{\{u\leq v\}}d\mu.
\eea
\end{lemma}

To establish this lemma,
recall that if
$\mu_j,\mu$ are non-negative Borel measures on a compact topological measure
space, with uniformly bounded total measures, and $d\mu_j\to d\mu$ weakly,
then for any open subset ${\cal O}$ and any compact subset $K$, we have
(see e.g. \cite{EG})
\bea
\label{weaklimits}
\int_{\cal O}d\mu\leq {\rm liminf}_{j\to\infty}\int_{\cal O}d\mu_j,
\qquad
{\rm limsup}_{j\to\infty}\int_{K}d\mu_k
\leq
\int_{K}d\mu.
\eea
If the measures $d\mu_j$ are uniformly bounded in capacity,
then the first inequality extends to all sets $E$ which are ``quasi-open'', in the sense that
for any $\e>0$, there exist an open set ${\cal O}$,
and Borel sets $G_\e\subset M'$ and $G_{\e}'\subset M'$ with capacities
less than $\e$ so that
\bea
E\subset {\cal O}\cup G_\e,\qquad
{\cal O}\subset E\cup G_\e'
\eea
Similarly,
the second inequality extends to all sets $E$ which are ``quasi-compact'', in the sense that
for any $\e>0$, there exist a compact set $K$,
and Borel sets $G_\e\subset M'$ and $G_{\e}'\subset M'$ with capacities
less than $\e$ so that $E\subset K\cup G_\e,\
K\subset E\cup G_\e'$.

\medskip
\noindent
{\it Proof of Lemma \ref{BT82-1}:}
We take limits in (\ref{comparisonineq1})
successively as $j\to\infty$ and then as $k\to\infty$.

\smallskip
First, consider the left hand side of (\ref{comparisonineq1}).
Both integrand and domain of integration depend on $j$,
so we change first to a domain of integration
independent of $j$ by writing
\bea
\int_{\{u_k<v_j\}}d\nu_j
\geq
\int_{\{u_k<v\}}d\nu_j
\eea
since $v_j\geq v$. Now the set $\{u_k<v\}$ is not necessarily open,
but it is quasi-open in the sense defined above. Indeed,
by the quasi-continuity of $v$, for each $\e>0$,
$v=V_\e$ for a function
$V_\e$ continuous on $M$, outside a set $G_\e$ of capacity less than $\e$.
Thus $\{u_k<V_\e\}\subset\{u_k<v\}\cup G_\e$ and
$\{u_k<v\}\subset\{u_k<V_\e\}\cup G_\e$,
and $\{u_k<V_\e\}$ is open.
Applying the inequality (\ref{weaklimits}) for quasi-open sets, we get
\bea
{\rm liminf}_{j\to\infty}\int_{\{u_k<v_j\}}d\nu_j
\geq
\int_{\{u_k<v\}}d\nu.
\eea
Next, the limit as $j\to\infty$ of the right hand side of (\ref{comparisonineq1})
can be bounded in a straightforward way by
\bea
{\rm lim}_{j\to\infty}
\int_{\{u_k<v_j\}}d\mu_k\leq \int_{\{u_k\leq v\}}d\mu_k.
\eea
Altogether, for each $k$,
the limit as $j\to\infty$ of the inequality (\ref{comparisonineq1}) produces
\bea
\label{intermediate}
\int_{\{u_k<v\}}d\nu
\leq
\int_{\{u_k\leq v\}}d\mu_k.
\eea

\medskip
The second step is to take the limit
of (\ref{intermediate}) as $k\to+\infty$.
The left hand side gives
\bea
{\rm lim}_{k\to\infty}\int_{\{u_k<v\}}d\nu
=
\int_{\{u<v\}}d\nu.
\eea
For the right hand side, where integrand and domain of integration both depend
on $k$, we argue in complete analogy with the preceding case
and begin by writing write
\bea
\int_{\{u_k\leq v\}}d\mu_k
\leq
\int_{\{u\leq v\}}d\mu_k,
\eea
since $u\leq u_k$.
The set $\{u\leq v\}$ is quasi-compact, since for each $\e$,
by the quasi-continuity of $u$, we can write
$u=U_\e$ outside a set of capacity less than $\e$,
with $U_\e$ a continuous function on $M$.
The weak convergence $d\mu_k\to d\mu$ implies, by (\ref{weaklimits}),
\bea
{\rm limsup}_{k\to\infty}
\int_{\{u_k\leq v\}}d\mu_k
\leq
{\rm limsup}_{k\to\infty}
\int_{\{u\leq v\}}d\mu_k
\leq
\int_{\{u\leq v\}}d\mu.
\eea
The lemma is proved.

\bigskip

We would like to apply Lemma \ref{BT82-1} to our context.
Let $\Omega_0$ be a smooth, non-negative closed $(1,1)$-form
on $\bar M$, and let $u,v$ be $\Omega_0$-plurisubharmonic and bounded
on the K\"ahler manifold $(M,\Omega)$. By the theorem
of Blocki and Kolodziej \cite{BK}, for any open subset $M'$
of $M$ with compact closure, there exists a decreasing
sequence $\e_j\downarrow 0$, and sequences $u_j$, $v_j$
of smooth functions with
\bea
u_j \downarrow u,
\quad
v_j\downarrow v
\eea
in an open neighborhood of $M'$,
and $u_j$ and $v_j$ are $(\Omega_0+\e_j\Omega)$-plurisubharmonic.

\begin{lemma}
\label{BT82-2}
Let $u,v\in L^\infty$ be $\Omega_0$-plurisubharmonic functions satisfying
\bea
\label{3delta}
{\rm liminf}_{z\to\partial M}
(u(z)-v(z))\geq 3\delta
\eea
for some fixed constant $\delta>0$.
Let $M'$ be any open subset of $M$ with compact closure,
with $u-v> 2\delta$ in a neighborhood $K$ of $\pl M'$.
Let $u_j$, $v_j$ be the decreasing sequences
approximating $u$ and $v$ as given by the theorem of Blocki and
Kolodziej, and let
\bea
&&
d\mu_k=\Omega_0+\e_k\Omega+{i\over 2}\ddb u_k,
\quad
d\nu_j=\Omega_0+\e_j\Omega+{i\over 2}\ddb v_j,
\nonumber\\
&&
d\mu=\Omega_0+{i\over 2}\ddb u, \quad
d\nu=\Omega_0+{i\over 2}\ddb v.
\eea
Then all five conditions (a-e) of Lemma \ref{BT82-1}
are satisfied.
\end{lemma}

\noindent
{\it Proof of Lemma \ref{BT82-2}:}
The condition (a) follows directly from the $\Omega_0$-plurisubharmonicity of $u,v$,
and the Bedford-Taylor Theorem on the quasi-continuity of plurisubharmonic
functions on ${\bf C}^m$, applied to each coordinate chart of $M$.

\smallskip
To prove (b), fix an index $k$.
For each point $z_0\in K$,
choose $j_{z_0}$ so that
$v_{j_0}(z_0)<u_k(z_0)-\delta$,
which is possible, since $v_j(z_0)$ converges to $v(z_0)$
and $v<u-2\delta\leq u_k-2\delta$.
By the upper-semicontinuity of the function $v_{j_0}-u_k$,
it follows that there is a neighborhood ${\cal O}_{z_0}$
with $v_{j_0}<u_k-\delta$ on ${\cal O}_{z_0}$.
Let $\cup_{\al=1}^N{\cal O}_{z_\al}$
be a finite cover of the compact set $K$,
and let $J_k={\rm max}_{1\leq\al\leq N}j_\al$.
Then for any $j\geq J_k$ and any $z\in K$,
pick $z_\al$ with $z\in {\cal O}_{z_\al}$. Since $v_j$
is a decreasing sequence, we have
\bea
v_j(z)\leq v_{j_\al}(z)
<u_k(z)-\delta,
\eea
which is the desired statement.

\smallskip
For (c), it suffices to establish the weak convergence on each compact subset of $M$.
Covering the compact set by a finite number of coordinate charts
$U_\alpha$, it suffices to establish the weak
convergence on each chart $U_\alpha$. We may assume that on $U_\alpha$,
\bea
\Omega_0={i\over 2}\ddb f_{0,\al},
\qquad
\Omega={i\over 2}\ddb f_{\al},
\eea
where $f_{\al}$ may be assumed $>0$ by adding a suitable large constant.
Then
\bea
0\leq \Omega_0+\e_j\Omega+{i\over 2}\ddb u_j
=
{i\over 2}\ddb(f_{0,\al}+\e_j f_\al+u_j)
\eea
Thus the functions $f_{0,\al}+\e_jf_{\al}+u_j$ are plurisubharmonic
and decreasing to $f_{0,\al}+u$. The Bedford-Taylor monotonicity theorem
implies the weak convergence on $M_\al$,
\bea
(\Omega_0+\e_j\Omega+{i\over 2}\ddb u_j)^m\to (\Omega_0+{i\over 2}\ddb u)^m.
\eea
The case of $v_j$ is similar, so this establishes (c).

\smallskip
The statement (d) is a consequence of the fact that $u_k$, $u$, $v_j$, and $v$
can be assumed to be all uniformly bounded in absolute value by the same constant $C$.
By the definition of capacity, it follows for example that
for each $E_\al\subset U_\al$, $U_\al$ coordinate chart, we have
\bea
\int_{E_\al} d\mu \leq \|f_{0,\al}+u\|_{L^\infty}^m\,c(E_\al,U_\al).
\eea

\smallskip
Finally, the statement (e) follows by applying the smooth version
Lemma \ref{BT76}, to the level sets $\{u_k+\lambda<v_j\}$,
which have compact closure in $M'$ and smooth boundary for generic $\lambda>0$.
Letting $\lambda\downarrow 0$ gives the desired inequality. Q.E.D.

\medskip
\noindent
{\it Proof of Theorem 4}:
If we replace $u$ by $u+3\delta$ with $\delta>0$,
then the condition (\ref{3delta})
is satisfied. Choosing $M'$ as in Lemma \ref{BT82-2},
we can apply Lemma \ref{BT82-1}, and obtain
\bea
\int_{\{u+3\delta<v\}}(\Omega_0+{i\over 2}\ddb v)^m
\leq
\int_{\{u+3\delta\leq v\}}(\Omega_0+{i\over 2}\ddb u)^m.
\eea
The theorem follows by letting $\delta\downarrow 0$.

\subsection{A uniqueness theorem for completely degenerate complex
Monge-Amp\`ere equations}

The comparison theorem implies the following uniqueness theorem,
for $\Omega_0$-plurisubharmonic solutions of a completely degenerate
Monge-Amp\`ere equations, where the form $\Omega_0$ is allowed to
be degenerate along an analytic subvariety:

\begin{theorem}
\label{uniqueness}
Let $(M,\Omega)$ be a K\"ahler manifold with smooth boundary
$\partial M$ and dimension $m$, and let $u,v\in L^\infty$
be $\Omega_0$-plurisubharmonic functions satisfying
\bea
(\Omega_0+{i\over 2}\ddb u)^m=(\Omega_0+{i\over 2}\ddb v)^m=0,
\qquad {\rm limsup}_{z\to\partial M}(u(z)-v(z))=0.
\eea
If $\Omega_0$ is $\geq 0$ everywhere, and $>0$ away from
an analytic subvariety of strictly
positive codimension which does not intersect $\partial M$,
then we must have $u=v$ on $M$.
\end{theorem}

\noindent
{\it Proof:} By adding the same large constant to
both $u$ and $v$, we may assume that $u,v>0$.
We argue by contradiction, and thus begin by
assuming that $S=\{u<v\}\not=\emptyset$.
Since $u,v$ are $\Omega_0$-plurisubharmonic, the set
$S$
must have strictly positive measure (it suffices to work in local
coordinates, and apply the corresponding well-known
property of plurisubharmonic functions on ${\bf C}^m$).
Furthermore, since we can write
\bea
S=\cup_{\e>0}\{u<(1-\e)v\}\equiv \cup_{\e>0}S_\e,
\eea
it follows that $S_\e$ must have strictly positive
measure for some $\e>0$. Fix one such value of $\e$.
Since $u\geq v\geq (1-\e)v$ on $\partial M$, we may apply
the comparison principle for K\"ahler manifolds and obtain
\bea
0\geq \int_{S_\e}(\Omega_0+{i\over 2}\ddb u)^m
&\geq&
\int_{S_\e}(\Omega_0+(1-\e){i\over 2}\ddb v)^m
\nonumber\\
&=&
\int_{S_e}\{(1-\e)(\Omega_0+{i\over 2}\ddb v)+\e\Omega_0\}^m
\nonumber\\
&\geq&
\e^m\int_{S_e}\Omega_0^m
\eea
since the form $\Omega_0+{i\over 2}\ddb v$ is non-negative. Let now $V_\delta$ be
the complement of a neighborhood of the divisor $D$, with $\Omega_0^m\geq \delta
\Omega^m$ for each $\delta>0$ small enough. Clearly for each $\delta>0$
\bea
\int_{S_\e}\Omega_0^m
\geq
\int_{S_\e\cap V_\delta}\Omega_0^m
\geq
\delta\int_{S_\e\cap V_\delta}\Omega^m.
\eea
Since $M\setminus D=\cup_{\delta>0}V_\delta$ and $D$ has measure $0$ with respect
to the volume form $\Omega^m$, we have
\bea
0<\int_{S_\e}\Omega^m
=
{\rm lim}_{\delta\to 0}\int_{S_\e\cap V_\delta}\Omega^m
\eea
which implies that $\int_{S_\e\cap V_\delta}\Omega^m>0$ for some $\delta>0$.
Altogether, we obtain a contradiction. Thus $\{u<v\}$ must be empty.
Interchanging the roles of $u$ and $v$ completes the proof of the theorem.

\section{Proof of Theorem 1}
\setcounter{equation}{0}

We can now prove Theorem 1. In Theorem \ref{MAOmega0theorem},
we have shown that the function $\hat\Phi$ is a bounded,
$\Omega_0$-plurisubharmonic solution of the
Dirichlet problem (\ref{MAOmega0}) on $\cX_D$.
On the other hand, in \cite{PS09} (Theorem 3),
it was shown that the same Dirichlet problem admits a bounded,
$\Omega_0$-plurisubharmomic solution which is $C^{1,\al}$
for any $0<\al<1$ on $\cX_D^\times$. By Theorem \ref{uniqueness},
it follows that the two solutions must coincide. Thus
$\hat\Phi$ is $C^{1,\al}$ on $\cX_D^\times$.
Since $\hat\Phi\ =\ \Phi-\Phi_1+\Psi_0$ and both $\Phi_1$ and $\Psi_0$
are smooth on $\cX_D^\times$, it follows that $\Phi$
is $C^{1,\al}$ on $\cX_D^\times=X\times D^\times$. Q.E.D.

\medskip
D.H. Phong

Department of Mathematics

Columbia University, New York, NY 10027

\bigskip

Jacob Sturm

Department of Mathematics

Rutgers University, Newark, NJ 07102

\end{document}